\documentclass[]{amsart}
\usepackage{amsmath,amssymb}
\newtheorem{theo}{\sc Theorem}

\def\C{\mathbb{C}}

\def\Nquer{\overline{N}}

\def\c{\mathfrak{c}}

\def\t{\mathfrak{t}}

\def\ende{ \square}
\def\Ende{~$\ende$\par}

\def\be{\begin{equation}}
\def\ee{\end{equation}}

\def\kzu1#1{\buildrel (#1) \over \longrightarrow}
\newcounter{ABS}

\def\ds{\displaystyle}
\def\ts{\textstyle}
\def\proof{\noindent {\it Proof.~}}
\author{\sc Norbert Steinmetz}
\title[]{Remark on meromorphic functions that share five pairs}
\address{Institut f\"ur Mathematik $\cdot$ TU Dortmund $\cdot$ D-44221 Dortmund $\cdot$ Germany}
 \email{stein@math.tu-dortmund.de}
 
\begin{document}

\maketitle

\begin{abstract}We determine all pairs $(f,g)$ of meromorphic functions  that share four pairs of values $(a_\nu,b_\nu)$, $1\le\nu\le 4$,
and a fifth pair $(a_5,b_5)$ under some mild additional condition.
  \end{abstract}

\section{Introduction}Meromorphic functions $f$ and $g$ are said to share the pair $(a,b)$ of complex numbers (including $\infty$), if
$f-a$ and $g-b$ ($1/f$ and $1/g$, if $a=\infty$ and $b=\infty$, respectively) have the same zeros.
Czubiak and Gundersen~\cite{CG} proved that meromorphic functions $f$ and $g$ that share {\it six} pairs $(a_\nu,b_\nu)$ are M\"obius transformations of each other, hence share all pairs $(a,L(a))$
for some M\"obius transformation $L$. On the other hand, the functions
\be\label{GEX}\hat f(z)=\ds\frac{e^z+1}{(e^z-1)^2}\quad{\rm and}\quad\hat g(z)=\frac{(e^z+1)^2}{8(e^2-1)}\ee
share the values $\infty,0,1,$ and $-\frac18$
with different multiplicities, and the pair $(-\frac12,\frac14)$ counting multiplicities. Thus
\be\label{GEXmod}f(z)=\frac{1}{\hat f(z)+\frac12}\quad{\rm and}\quad g(z)=\frac{1}{\hat g(z)-\frac14}\ee
are not M\"obius transformations of each other and share the pairs $(0,0)$, $(2,-4)$, $(\frac23,\frac43)$ and $(\frac83,-\frac83)$ with different multiplicities, and the value $\infty$ (the pair $(\infty,\infty)$) counting multiplicities. Moreover, $f$ and $g$ have common counting function of poles $\Nquer(r,\infty)=T(r)+S(r)$, where $T(r)$ and $S(r)$ denote the common Nevanlinna characteristic and remainder term of $f$ and $g$ (for notations and results of Nevanlinna theory the reader is referred to Hayman's monograph~\cite{hayman}), and $f$ and $g$ parametrise the algebraic curve
\be\label{GEXalg}4x^2+2xy+y^2-8x=0.\ee
Gundersen's example $\hat f$, $\hat g$ was the first to show that in Nevanlinna's {\it Four Value Theorem} \cite{nevanlinna2} one cannot dispense with the condition `counting multiplicities' for each of the four values. This is possible for one value (Gundersen~\cite{gundersen1}) and also for two of the values (Gundersen~\cite{gundersen2}, Mues~\cite{mues}), while the case of three such values is still open. The state of art is outlined in \cite{NSt4}.
Gundersen's example also has another characterisation due to Reinders~\cite{reinders1,reinders2}:
{\it If $f$ and $g$ share four values $a_\nu$, and if $f^{-1}(a)\subset g^{-1}(b)$ holds for some pair
$(a,b)$ $(a,b\ne a_\nu)$, then either $f$ and $g$ are M\"obius transformations of each other or else
$f=T\circ\hat f\circ h$ and $g=T\circ \hat g\circ h$ holds for some M\"obius transformation $T$ and some non-constant entire function $h$.}

In \cite{GGG}, Gundersen considered functions $f$ and $g$ that share five pairs and are not M\"obius transformations of each other. He proved several sharp inequalities for the corresponding Nevanlinna functions, including
$T(r,f)=T(r,g)+S(r)$ and
\be\label{Gund13}\Nquer(r;a_\nu,b_\nu)\ge\ts\frac13 T(r)+S(r),\ee
where $\Nquer(r;a_\nu,b_\nu)$ denotes the counting functions of common $(a_\nu,b_\nu)$-points of $(f,g)$,
not counting multiplicities, and $T(r)$ and $S(r)$ denote the common Nevanlinna characteristic and remainder term, respectively.

\section{Main result}
The aim of this paper is to prove

\begin{theo}Suppose that meromorphic functions $f$ and $g$ share four pairs $(a_\nu,b_\nu)$, and a fifth pair $(a_5,b_5)$ counting multiplicities and such that
\be\label{GenAss}m\big(r,1/(f-a_5)\big)+m\big(r,1/(g-b_5)\big)=S(r)\ee
holds. Then either $f$ and $g$ are M\"obius transformations of each other or else
$f=T\circ\hat f\circ h$ and $g=S\circ \hat g\circ h$ holds for suitably chosen M\"obius transformations $S$ and $T$ and some non-constant entire function $h$.\end{theo}

\proof We note that (\ref{GenAss}) is automatically fulfilled if $a_\nu=b_\nu$, $1\le\nu\le 4$. Three of the pairs $(a_\nu,b_\nu)$ may be prescribed. We will assume $(a_1,b_1)=(0,0)$, $(a_2,b_2)=(2,-4)$, and, in particular, $a_5=b_5=\infty$, to stay as close as possible with the modified example of Gundersen. Then  $f$ and $g$ have the same poles counting multiplicities, and such  that
\be\label{GeneralAssumption}m(r,f)+m(r,g)=S(r)\ee
holds. We also assume that $f$ and $g$ are not M\"obius transformations of each other. Similar to the approach in \cite{GGG} we consider
$$P(x,y,\c)=c_1x^2+c_2xy+c_3y^3+c_4x+c_5y.$$
Then there are at least two linear independent vectors $\c=(c_1,\ldots c_5)\in\C^5$ such that
\be\label{BED1}P(a_\nu,b_\nu,\c)=0\quad(1\le\nu\le 4)\ee
holds, that is, $P(z)=P(f(z),g(z),\c)$
vanishes whenever $f(z)=a_\nu$ and $g(z)=b_\nu$. If $P$ does not vanish identically, this yields
$$\sum_{\nu=1}^4 \Nquer(r;a_\nu,b_\nu)\le \Nquer(r,1/P)\le T(r,P)+O(1)\le 2T(r)+S(r);$$
for the last inequality the additional hypothesis (\ref{GeneralAssumption}) is used.
On the other hand it follows from the Second Main Theorem that
$$\sum_{\nu=1}^4 \Nquer(r;a_\nu,b_\nu)+\Nquer(r,\infty)\ge 3T(r)+S(r),$$
hence $T(r)\le \Nquer(r,\infty)+S(r).$
Thus, still assuming $P\not\equiv 0$, it follows that
$$\begin{array}{rcl}
\Nquer(r,1/P)&=&N(r,1/P)+S(r)=2T(r)+S(r)\cr
m(r,1/P)&=&S(r)\cr
\sum\limits_{\nu=1}^4\Nquer(r,a_\nu,b_\nu)&=&\Nquer(r,1/P)+S(r)\cr
T(r)&=&\Nquer(r,\infty)+S(r).\end{array}$$
In particular, the quotient $\chi(z)=P(z;\tilde\c)/P(z;\c)$
satisfies $T(r,\chi)=S(r).$
In other words, $f$ and $g$ parametrise the algebraic curve
\be\label{Algcurvechi}F(x,y;z)=\chi_1x^2+\chi_2yx+\chi_3y^2+\chi_4x+\chi_5y=0\quad(\chi_k=\chi c_k-\tilde c_k)\ee
over the field $\C(\chi)$. This is also true if $P(z;\c)$ or $P(z;\tilde\c)$ vanishes identically. It is obvious that $\chi_1\chi_3\not\equiv 0$, since otherwise
$g$ [resp. $f$] would be a M\"obius transformation or a rational function of $f$ [resp. $g$] of degree two over the field $\C(\chi)$. In the first case it would follow that $g$ is an ordinary M\"obius transformation of $f$, while in the second case we would obtain a contradiction: $T(r,g)=2T(r,f)+S(r)$.

The algebraic curve (\ref{Algcurvechi}) has the rational parametrisation (set $x=ty$)
$$x=\frac{p(z,t)}{s(z,t)}=-\frac{t(\chi_4t+\chi_5)}{\chi_1t^2+\chi_2t+\chi_3},
~y=\frac{q(z,t)}{s(z,t)}=-\frac{\chi_4t+\chi_5}{\chi_1t^2+\chi_2t+\chi_3}$$
with $t=x/y.$ In terms of $f$ and $g$ this yields
$$\begin{array}{rcl}
f(z)&=&\ds\frac{p(z,\t(z))}{s(z,\t(z))}=-\frac{\t(z)(\chi_4\t(z)+\chi_5)}{\chi_1\t(z)^2+\chi_2\t(z)+\chi_3}\cr
g(z)&=&\ds\frac{q(z,\t(z))}{s(z,\t(z))}=-\frac{\chi_4\t(z)+\chi_5}{\chi_1\t(z)^2+\chi_2\t(z)+\chi_3}\end{array}\quad{\rm with}\quad\t(z)=\frac{f(z)}{g(z)}.$$
Since by (\ref{Gund13}), $f$ and $g$ have `many' zeros, there are three possibilities to be discussed: The zeros correspond to the
\begin{itemize}
\item[a)] {\it poles of $\t$,} in which case $\chi_4\equiv 0$ and `almost all' zeros of $f$
are simple, while the zeros of $g$ have order two. Moreover, $\t$ has `almost no' zeros ($N(r,1/\t)=S(r)$).
\item[b)] {\it zeros of $\t$}, in which case $\chi_5\equiv 0$ and `almost all' zeros of $g$ are
simple, while the zeros of $f$ have order two. Moreover, $\t$ has `almost no' poles ($N(r,\t)=S(r)$).
\item[c)] {\it zeros of $\chi_4(z)\t(z)+\chi_5(z)$ with $\chi_4\chi_5\not\equiv 0$.} Then
`almost all' zeros of $f$ and $g$ are simple, and $\t$ has `almost no' zeros and poles ($N(r,1/\t)+N(r,\t)=S(r)$).
\end{itemize}
Taking all pairs $(a_\nu,b_\nu)$ ($1\le\nu\le 4$),into account, the following holds: for every $\nu$ there exist $\phi_\nu,\psi_\nu,\alpha_\nu,\beta_\nu,\tilde\beta_\nu\in\C(\chi)$ such that
$p(z,t)-a_\nu s(z,t)=\phi_\nu (t-\alpha_\nu )(t-\beta_\nu )$ and $q(z,t)-b_\nu s(z,t)=\psi_\nu (t-\alpha_\nu )(t-\tilde\beta_\nu )$, respectively;  
occasionally the factor $(t-\beta_\nu)$ and $(t-\tilde\beta_\nu)$
corresponding to $\beta_\nu\equiv\infty$ and $\tilde\beta_\nu\equiv\infty$, respectively, might be missing.
The functions(\footnote{At first glance one would expect that $\alpha_\nu,\beta_\nu,\tilde\beta_\nu$ are algebraic over $\C(\chi)$. But this is not the case, since analytic continuation which permutes $\alpha_\nu$ and $\beta_\nu$ would also permute $\alpha_\nu$ and $\tilde\beta_\nu$, in contrast to $\beta_\nu\not\equiv\tilde\beta_\nu$.}) $\alpha_\nu$ are mutually distinct, and the same is true for $\beta_\nu$ and also $\tilde\beta_\nu$. It is also obvious that $\beta_\nu\not\equiv\tilde\beta_\nu$,
and that both functions are exceptional for $\t$, except when one of them coincides with $\alpha_\nu$. Since $\t$ has at most two exceptional functions, we obtain the following picture:

For $\nu=1$ and $\nu=4$, say, we have $\beta_\nu\equiv\alpha_\nu$, that is, the pairs $(a_\nu,b_\nu)$, are attained by $(f,g)$ in a $(2:1)$ manner, while for $\nu=2$ and $\nu=3$ this happens the other way $(1:2)$. This means that, in addition to (\ref{Algcurvechi}), that also
\be\label{BED2}F_y(a_\nu,b_\nu;z)\equiv 0\quad(\nu=1,4)\quad{\rm and}\quad F_x(a_\nu,b_\nu;z)\equiv 0\quad(\nu=2,3)\ee
holds. To stay close with the modified example of Gundersen we assume $\chi_3\equiv 1$
(this is possible since $\chi_3\not\equiv 0$ is already known). From  (\ref{BED2}), that is
$$\chi_5=4\chi_1-4\chi_2+\chi_4=2\chi_1a_3+\chi_3b_3+\chi_4=\chi_2a_4+2b_4\equiv 0,$$
one can compute the coefficients $\chi_k$ in terms of $a_3$, $b_3$, $a_4$, $b_4$, namely
\be\label{KOEFF}\chi_1=\frac{b_4(b_3+4)}{a_4(a_3-2)},~\chi_2=-\frac{2b_4}{a_4},~\chi_3=1,~
\chi_4=\frac{2b_4(2b_3+4a_3)}{a_4(2-a_3)}.\ee
In particular, the functions $\chi_k$ are constant, and $f$ and $g$ are rational functions (now over $\C$) of the meromorphic function $\t=f/g$.
Having determined the coefficients (\ref{KOEFF}) we now use (\ref{Algcurvechi}) to express $b_3$ and $b_4$ in terms of $a_3$ and $a_4$.
The solutions to $F(a_\nu,b_\nu;z)=0$ for $\nu=2,4$ are given by(\footnote{We note that {\sf maple} was not able to determine all solutions to the system (\ref{BED1}), but note also that the coefficients $\chi_k$ are functions of the $a_\nu$, $b_\nu$; to get an impression: one has to solve
$$(b_3+4)b_4+(8-4a_3)a_4=a_3^2b_4(b_3+4)+b_3^2(2-a_3)=a_4(b_3+4)+a_3(8-b_4)+2b_4+4b_3=0$$
for $b_3,b_4$.})

$\bullet$  $b_4 = -2a_4, b_3 = -2a_3$, while $F(a_3,b_3;z)=0$ is automatically fulfilled; this leads to $g=-2f$.

$\bullet$ $b_4=2a_4-8$ and $b_3=\ds\frac{2(8-4a_4+a_4a_3)}{a_4-4}$.
Since, however, $\ds F(a_3,b_3;z)=32\frac{(a_3-2)(a_4-2)(a_3-a_4+2)}{(a_4-4)^2}$
also has to vanish, we just have to discuss the sub-case $a_3=a_4-2$, since $a_3=2$ and also $a_4=2$ would contradict $a_2=2$.
Thus $a_3=a_4-2$, $b_3=2a_4-4$ and $b_4=2a_4-8$, and $\ts(a_1,a_2,a_3,a_4)$ and $(b_1,b_2,b_3,b_4)$
have the same cross-ratio $\ds\frac{a_1-a_3}{a_2-a_3}:\frac{a_1-a_4}{a_2-a_4}=\frac{b_1-b_3}{b_2-b_3}:\frac{b_1-b_4}{b_2-b_4}=\frac{(a_4-2)^2}{a_4(a_4-4)}$.
In other words, there exists some  M\"obius transformation $L$ such that $f$ and $L\circ g$ share four {\it values}
$a_1,a_2,a_3,a_4$ and the pair $(\infty,L(\infty))$. By Reinders' characterisation this implies $f=T\circ \hat f\circ h$ and $g=S\circ \hat g\circ h$, where $S$ and $T$ are suitably chosen M\"obius transformations, and $h$ is some non-constant entire function. \Ende

\medskip{\sc Final remark.} It remains open whether or not--and how--the hypothesis (\ref{GeneralAssumption}) may be relaxed. Is it sufficient to assume that the pair $(a_5,b_5)$ is shared  `counting multiplicities' by $f$ and $g$?  Is it even true that functions sharing five pairs are either M\"obius transformations of each other or else have the form
$f=T\circ\hat f\circ h$ and $g=S\circ \hat g\circ h$?

\end{document}